\documentclass[12pt]{amsart}
\usepackage{epsf,amsmath,amscd,amsthm,amsfonts,latexsym,epsfig,graphicx}

\theoremstyle{plain}
\newtheorem{theorem}                 {Theorem}      [section]
\newtheorem{proposition}  [theorem]  {Proposition}
\newtheorem{corollary}    [theorem]  {Corollary}
\newtheorem{lemma}        [theorem]  {Lemma}

\theoremstyle{definition}

\newtheorem{remark}       [theorem]  {Remark}

\numberwithin{equation}{section}

\def \r{\mbox{${\mathbb R}$}}
\def \s{\mbox{${\mathbb S}$}}
\def \e{\mbox{${\mathbb E}$}}
\def \h{\mbox{${\mathbb H}$}}

\DeclareMathOperator{\trace}{trace}

\begin{document}

\begin{abstract}
We consider the biharmonicity condition for maps between
Riemannian manifolds (see \cite{PBDK}), and study the non-geodesic
biharmonic curves in the Heisenberg group $\h_3$. First we prove
that all of them are helices, and then we obtain explicitly their
parametric equations.

\end{abstract}

\title[Explicit formulas for non-geodesic biharmonic curves...]{Explicit formulas for
 non-geodesic\\  Biharmonic curves of the Heisenberg group}


\author{R. Caddeo}
\address{Universit\`a degli Studi di Cagliari\\
Dipartimento di Matematica\\
Via Ospedale 72\\
09124 Cagliari\\
ITALIA}
\email{caddeo@unica.it}
\author{C. Oniciuc}
\address{Faculty of Mathematics\\"AL.I.~Cuza" University of Iasi\\
Bd. Carol I, no. 11\\
700506 Iasi\\
ROMANIA}
\email{oniciucc@uaic.ro}
\author{P. Piu}
\address{Universit\`a degli Studi di Cagliari\\
Dipartimento di Matematica\\
Via Ospedale 72\\
09124 Cagliari\\
ITALIA}
\email{piu@unica.it}


\maketitle

\section{Introduction}

By the definition given by J.~Eells and J.H.~Sampson in \cite{JEJS},
a map $\phi$ from a compact Riemannian
manifold $(M,g)$ to another Riemannian
manifold $(\overline{M},\overline{g})$ is harmonic
if it is a critical point of the energy
\[
E_1(\phi)=\frac{1}{2}\int_M
\vert d\phi\vert^2v_g.
\]
\vspace*{0.2 cm}\\
From the first variation formula it follows that $\phi$ is
harmonic if and only if its first tension field
$\tau_1(\phi)=\trace \nabla d\phi$ vanishes. The same authors
suggested a generalization of the notion of harmonicity: a map
$\phi$ is biharmonic if it is a critical point of the bienergy
functional
\[
E_2(\phi)=\frac{1}{2}\int_M\vert\tau_1(\phi)\vert^2v_g.
\]
\vspace*{0.2 cm}\\
The first variation formula for $E_2$, given by G.Y.~Jiang in
\cite{GYJ1} and \cite{GYJ2}, is reobtained in \cite{PBDK} and can be written as
\vspace*{0.2 cm}
$$
\frac{dE_2(\phi_t)}{dt}{\biggr \vert}_{t=0}=\int_M<\tau_2(\phi),V>v_g,
$$
\vspace*{0.2 cm}\\
where $\{\phi_t\}_{t}$ is a variation of $\phi$, $V$ is the
variational vector field along $\phi$, and
\vspace*{0.2 cm}
\begin{align}\label{eq:bihar_1}
\tau_2(\phi) &= - J(\tau_1(\phi)) \\
&= -\Delta\tau_1(\phi) +\trace R^{\overline{M}}(d\phi,\tau_1(\phi))d\phi.\notag
\end{align}
\vspace*{0.2 cm}\\
Here $J$ denotes the Jacobi operator and
 $\Delta\tau_1(\phi) = - \trace (\nabla^{\phi})^2 \tau_1(\phi)$, where
$\nabla^{\phi}$
represents the connection in $\phi^{-1} (T\overline{M})$.
Hence the condition of biharmonicity for $\phi$ is
\vspace*{0.2 cm}
\begin{equation}\label{tensione2}
\tau_2(\phi)=0.
\end{equation}
\vspace*{0.2 cm}\\
When $M$ is not compact, as a definition of biharmonicity one can
adopt equation (\ref{tensione2}) instead of the definition given
in terms of the integral formula of the bienergy.

Since any harmonic map is evidently
biharmonic, we are interested in non-harmonic biharmonic maps.

In general, the biharmonic equation is very complicated, and
therefore the problem of finding examples of non-harmonic
biharmonic maps between two Riemannian manifolds, or examples of
non-harmonic (non-minimal) biharmonic submanifolds of a given
Riemannian manifold, is difficult to solve. Still there are some
results that now we mention.
\begin{itemize}
\item
Biharmonic curves on a surface in $\r^3$ have been considered in
\cite{RCSMPP1}.
\item In \cite{RCSMCO1} the authors gave the complete classification of
non-minimal biharmonic submanifolds of $\s^3$. These are certain
circles, spherical helices and parallel spheres. The general case
of $\s^n$, for $n>3$, is more interesting and has been treated in
\cite{RCSMCO2}.
\item More recently, in \cite{PBDK}, it has been shown how to construct
non-harmonic
biharmonic maps from $M$ to $\overline{M}$ by deforming conformally
the metric of $M$.
\item When $M$ is a submanifold of the Euclidean space $\r^n$, the biharmonicity
condition seems to be very restrictive. In fact, B.Y. Chen
conjectured in \cite{BYCSI} that any biharmonic submanifold of the
Euclidean space $\r^n$ is minimal. This was proved to be true in
the case of submanifolds in $\r^3$ (in \cite{BYCSI}) and in other
special cases (see \cite{ID}, \cite{THTV} and \cite{BYC}).
\item The same result was proved in \cite{RCSMCO2} for any $3$-dimensional
Riemannian
manifold
$\overline{M}^{3}(-1)$ with constant negative sectional curvature
$-1$.
\item Other results about non-existence of non-harmonic biharmonic maps, when
the sectional curvature or the Ricci curvature of $\overline{M}$
is non-positive, can be found in \cite{GYJ1}, \cite{CO1} and
\cite{CO2}.
\end{itemize}
It seems then natural, as the next step, to consider  biharmonic
submanifolds in a $3$-manifold $\overline{M}$ with non-constant
sectional curvature. We choose as $\overline{M}$ the Heisenberg
group $\h_3$. This nilpotent Lie group is not-symmetric;
nevertheless it has many symmetries, in the sense that its
isometry group has dimension four, the biggest possible for a
3-space of non-constant sectional curvature. It is probably in
virtue of these symmetries that some main problem turns out to be
easier than expected. As, for instance, the problem of finding the
equations of geodesics (see \cite{OK} and \cite{ThH}), and that of
determining all surfaces with constant mean or Gaussian curvature
invariant with respect to some subgroup of the isometry group (see
\cite{PT}, \cite{RCPPAR1}, \cite{RCPPAR2} and \cite{CFFMRP}). This
manifold has many interesting properties. We mention, for example,
that there are no totally umbilical surfaces and therefore there
are no totally geodesic surfaces in $\h_3$ (see \cite{AS}).

In this paper we first write down the conditions that any non-harmonic
(non-geodesic)
biharmonic curve in
$\h_3$
must satisfy. Then we prove that the non-geodesic biharmonic curves
 in $\h_3$ are helices. A similar fact occurs
in $\s^3$. Finally we deduce the explicit parametric equations of the
non-geodesic
biharmonic curves in $\h_3$.
\vspace{0.3 cm}

{\it Notation.} We shall work in $C^{\infty}$ category, i.e.
manifolds, metrics, connections, maps will be assumed to be
smooth.

\section{Riemannian structure of $\h_3$}

The Heisenberg group $\h_3$ can be seen as the Euclidean space
$\r^3$ endowed with the multiplication
\vspace*{0.2 cm}
\[
(\widetilde{x},\widetilde{y},\widetilde{z})(x,y,z)=
(\widetilde{x}+x,\widetilde{y}+y,\widetilde{z}+z+
\frac{1}{2}\widetilde{x}y-
\frac{1}{2}\widetilde{y}x)
\]
\vspace*{0.2 cm}\\
and with the Riemannian
metric $g$ given by
\vspace*{0.2 cm}
\begin{equation}\label{metricaheis}
g=dx^2+dy^2+(dz+\frac{y}{2}dx-\frac{x}{2}dy)^2.
\end{equation}
\vspace*{0.2 cm}\\
The metric $g$ is invariant with respect to the left-translations
corresponding to that
multiplication. This metric is isometric to the other, also quite standard,
which is
left-invariant with respect to the composition arising from the
multiplication of the
$3\times3$ Heisenberg matrices.

At each point the metric $g$ has an axial symmetry: the 4-dimensional
group of its isometries
contains the group of rotations around the $z$ axis (in the classical
terminology (see \cite{LB}) a
space with such a property is called {\it systatic}) .\\

First of all we shall determine the Levi-Civita connection $\nabla$ of the
metric
$g$ with
respect to the left-invariant orthonormal basis
\vspace*{0.2 cm}
\[
e_1=\frac{\partial}{\partial x}-\frac{1}{2}y
\frac{\partial}{\partial z},
\quad
e_2=\frac{\partial}{\partial y}+\frac{1}{2}x
\frac{\partial}{\partial z},
\quad
e_3=\frac{\partial}{\partial z},
\]
\vspace*{0.2 cm}\\
which is dual to the coframe
\vspace*{0.2 cm}\\
\[
\theta^1=dx, \quad \theta^2=dy, \quad \theta^3=dz+\frac{y}{2}dx-
\frac{x}{2}dy.
\]
\vspace*{0.2 cm}\\
We obtain
\begin{equation}
\left\{
\begin{array}{l}
\nabla_{e_1}e_1=0, \quad \nabla_{e_1}e_2=\frac{1}{2}e_3, \quad
\nabla_{e_1}e_3=-\frac{1}{2}e_2,
\\ \mbox{} \\
\nabla_{e_2}e_1=-\frac{1}{2}e_3, \quad \nabla_{e_2}e_2=0, \quad
\nabla_{e_2}e_3=\frac{1}{2}e_1,
\\ \mbox{} \\
\nabla_{e_3}e_1=-\frac{1}{2}e_2, \quad \nabla_{e_3}e_2=\frac{1}{2}e_1,
\quad
\nabla_{e_3}e_3=0.
\end{array}
\right.
\label{eq:conn_1}
\end{equation}
\vspace*{0.2 cm}\\
Also, we have the well known Heisenberg bracket relations
\vspace*{0.2 cm}
\[
[e_1,e_2]=e_3, \quad [e_3,e_1]=[e_2,e_3]=0.
\]
\vspace*{0.2 cm}\\
We shall adopt the following notation and sign convention. The curvature
operator is given by
\vspace*{0.2 cm}\\
\[
R(X,Y)Z=-\nabla_X\nabla_YZ+\nabla_Y\nabla_XZ+\nabla_{[X,Y]}Z,
\]
\vspace*{0.2 cm}\\
while  the Riemann-Christoffel tensor field and the Ricci tensor field
are given by
\vspace*{0.2 cm}
\[
R(X,Y,Z,W)=g(R(X,Y)Z,W), \quad
\rho(X,Y)=\trace(Z\to R(X,Z)Y),
\]
\vspace*{0.2 cm}\\
where $X,Y,Z,W$ are smooth vector fields on $\h_3$.
Moreover we put, as in \cite{PPAS} or \cite{AS},
\vspace*{0.2 cm}
\[
R_{abc}=R(e_a,e_b)e_c, \quad R_{abcd}=R(e_a,e_b,e_c,e_d), \quad
\rho_{ab}=\rho(e_a,e_b),
\]
\vspace*{0.2 cm}\\
where the indices $a,b,c,d$ take the values $1,2,3$. The non vanishing
components
of the above tensor fields are
\vspace*{0.2 cm}
\begin{equation}
\left\{
\begin{array}{l}
R_{121}=-\frac{3}{4}e_2, \quad R_{131}=\frac{1}{4}e_3,
\quad R_{122}=\frac{3}{4}e_1,
\quad R_{232}=\frac{1}{4}e_3, \\ \mbox{} \\
R_{133}=-\frac{1}{4}e_1, \quad R_{233}=-\frac{1}{4}e_2,
\end{array}
\right.
\label{eq:curb_1}
\end{equation}
\begin{equation}
\label{eq:curb_2}
R_{1212}=-\frac{3}{4}, \quad R_{1313}=R_{2323}=\frac{1}{4},
\end{equation}
\begin{equation}
\label{eq:Ricci_1}
\rho_{11}=\rho_{22}=-\frac{1}{2}, \quad \rho_{33}=\frac{1}{2},
\end{equation}
\vspace*{0.2 cm}\\
respectively, and those obtained from these by means of the symmetries of $R$.
Thus the curvatures of $\h_3$ have both positive and negative components.
\vspace*{0.4 cm}
\section{Biharmonic curves in $\h_3$}

To study the biharmonic curves in $\h_3$, we shall use their
Frenet vector fields and equations. Let $\gamma:I\to \h_3$ be a
differentiable curve parametrized by arc length and let
$\{T,N,B\}$ be the orthonormal frame field tangent to $\h_3$ along
$\gamma$ and defined as follows: by $T$ we denote the unit  vector
field $\gamma'$ tangent to $\gamma$, by $N$ the unit vector field
in the direction of $\nabla_{T}T$ normal to $\gamma$, and we
choose $B$ so that $\{T,N,B\}$ is a positive oriented orthonormal
basis. Then we have the following Frenet equations \vspace*{0.2
cm}
\begin{equation}\label{eq:frenet}
\begin{array}{lcccc}
\nabla_{T}T&=&&kN&\\
\nabla_{T}N&=&-kT&& -\tau B\\
\nabla_{T}B&=&&\tau N&
\end{array},
\end{equation}
\vspace*{0.2 cm}\\
where $k=|\tau_1(\gamma)|=|\nabla_{T}T|$ is the geodesic curvature of
$\gamma$ and $\tau$ its geodesic torsion. By making use of equations
\eqref{eq:frenet} and of expression \eqref{eq:curb_1} of the
curvature tensor field, we obtain from \eqref{eq:bihar_1} the
biharmonic equation for $\gamma$:
\vspace*{0.2 cm}
\begin{eqnarray*}
\tau_{2}(\gamma)&=&\nabla_{T}^{3}T+R(T,kN)T \\
&=&(-3k'k)T+(k''-k^3-k\tau^2+\frac{k}{4}-kB_3^2)N
+(-2k'\tau-k\tau'+kN_3B_3)B\\
&=&0\,,
\end{eqnarray*}
\vspace*{0.2 cm}\\
where  $T = T_1e_1 + T_2e_2 + T_3e_3$, $N = N_1e_1 + N_2e_2 + N_3e_3$,
and $B = T \times N = B_1e_1 + B_2e_2 + B_3e_3$.
Thus we have

\begin{theorem}
\label{eq:theo_1}
Let $\gamma:I\to \h_3$ be a differentiable curve parametrized by arc
length.
Then $\gamma$ is a non-geodesic biharmonic curve if and only if
\begin{equation}\label{3.3}
\left\{
\begin{array}{l}
k=constant\neq 0 ;\\
k^2+\tau^2=\frac{1}{4}-B_3^2 ;\\
\tau'=N_3B_3 .
\end{array}
\right.
\end{equation}
\end{theorem}

\begin{remark} By analogy with  curves in $\r^3$, also following \cite{GL},
we keep
the name {\it helix}
for a curve in a Riemannian manifold having constant both geodesic curvature
and geodesic torsion. Now, for any helix in $\h_3$,
the system \eqref{3.3} becomes
\vspace*{0.2 cm}
\begin{equation}
\left\{
\begin{array}{l}
k^2+\tau^2=\frac{1}{4}-B_3^2 ;\\
N_3B_3 = 0 ,
\end{array}
\right.
\end{equation}
\vspace*{0.2 cm}\\
and therefore, in this case,  $B_3$ must be constant, too.
\end{remark}
Thus biharmonic helices satisfy
\begin{equation}\label{elic-biarm}
\left\{
\begin{array}{l}
B_3 = constant ;\\
k^2+\tau^2=\frac{1}{4}-B_3^2 ;\\
N_3 B_3 = 0 .
\end{array}
\right.
\end{equation}
\vspace*{0.2 cm}\\
We shall come back to biharmonic helices in the next section,
after showing that all the biharmonic curves in $\h_3$ are helices.
First we prove that for a biharmonic curve in $\h_3$ the geodesic torsion
must be constant.
\begin{proposition}\label{gener-prop}
Let $\gamma:I\to \h_3$ be a non-geodesic curve parametrized by arc
length. If $k$ is constant and $N_3B_3\neq 0$, then $\gamma$ is not
biharmonic.
\end{proposition}

\begin{proof}
From
\begin{eqnarray*}
\nabla_TT&=&(T_1'+T_2T_3)e_1+(T_2'-T_1T_3)e_2+T_3'e_3 \\
&=&kN,
\end{eqnarray*}
we obtain $T_3'=kN_3$; then, if we put $T_3(s)=kF(s)$ and $f(s) =F'(s)$,
we get
$N_3(s)=f(s)$. Hence we write
\vspace*{0.2 cm}
\[
T=\sqrt{1-k^2F^2}\cos\beta(s) e_1+\sqrt{1-k^2F^2}\sin\beta(s) e_2+kF(s)
e_3
\]
\vspace*{0.2 cm}\\
and then we use the first Frenet equation, that gives
\vspace*{0.2 cm}
\begin{eqnarray*}
\nabla_TT&=&\Big( \sqrt{1-k^2F^2}\big(kF-\beta'\big)\sin\beta -
\frac{k^2Ff}{\sqrt{1-k^2F^2}}\cos\beta \Big)e_1 \\
&&-\Big( \sqrt{1-k^2F^2}\big(kF-\beta'\big)\cos\beta +
\frac{k^2Ff}{\sqrt{1-k^2F^2}}\sin\beta \Big)e_2+kfe_3 \\
&=&kN.
\end{eqnarray*}
\vspace*{0.2 cm}\\
Since $k^2=\vert\nabla_TT\vert^2$, we have
$$
kF-\beta'=\pm k\frac{\sqrt{1-f^2-k^2F^2}}{1-k^2F^2}.
$$
\vspace*{0.2 cm}\\
Now we replace $kF-\beta'$ in the above expression of $\nabla_TT$, and we
obtain
\vspace*{0.2 cm}\\
\begin{eqnarray*}
N&=&\Big( \pm\frac{\sqrt{1-f^2-k^2F^2}}{\sqrt{1-k^2F^2}}\sin\beta
-\frac{kFf}{\sqrt{1-k^2F^2}}\cos\beta \Big)e_1 \\
&&+\Big( \mp\frac{\sqrt{1-f^2-k^2F^2}}{\sqrt{1-k^2F^2}}\cos\beta
-\frac{kFf}{\sqrt{1-k^2F^2}}\sin\beta \Big)e_2+fe_3.
\end{eqnarray*}
\vspace*{0.2 cm}\\
As $B=T\times N$, we have
$B_3=T_1N_2-N_1T_2=\mp\sqrt{1-f^2-k^2F^2}.$
Then the second Frenet equation gives
\vspace*{0.2 cm}
\begin{equation}
\label{eq:frenetequation1} <\nabla_TN,e_3>=<-kT-\tau
B,e_3>=-kT_3-\tau B_3.
\end{equation}
\vspace*{0.2 cm}\\
On the other hand we have
\vspace*{0.2 cm}
\begin{eqnarray}
\label{eq:frenetequation2}
<\nabla_TN,e_3>&=&<\nabla_T(N_1e_1+N_2e_2+N_3e_3),e_3>
\nonumber \\
&=&<\big(N_1'+\frac{1}{2}(T_2N_3+T_3N_2)\big)e_1+
\big(N_2'-\frac{1}{2}(T_1N_3+T_3N_1)\big)e_2
\nonumber \\
&&+\big(N_3'+\frac{1}{2}(T_1N_2-N_1T_2)\big)e_3,e_3>
\nonumber \\
&=&N_3'+\frac{1}{2}B_3.
\end{eqnarray}
\vspace*{0.2 cm}\\
By comparing ~\eqref{eq:frenetequation1} and ~\eqref{eq:frenetequation2}
we obtain
\begin{equation}
\label{eq:frenetequation3}
 N_3'+\frac{B_3}{2}=-kT_3-\tau B_3.
\end{equation}
\vspace*{0.2 cm}\\
Next we  replace $N_3=f$, $B_3=\mp\sqrt{1-f^2-k^2F^2}$ and $T_3=kF$ in
~\eqref{eq:frenetequation3}, and  we get
\vspace*{0.2 cm}
\begin{equation}
\label{eq:relation}
\tau = \pm \frac{f' + k^{2} F}{\sqrt{1-f^2-k^2F^2}} -\frac{1}{2}=
\frac{B_3'}{N_3}-\frac{1}{2}.
\end{equation}
\vspace*{0.2 cm}\\
Assume now that $\gamma$ is biharmonic. Then $\tau' = N_3B_3\neq0$
and we can write
\[
N_3=\frac{\tau'}{B_3} .
\]
 By substituting $N_3$ in ~\eqref{eq:relation} and then
by integrating  we get
\vspace*{0.2 cm}
\begin{equation}
\label{eq:int}
\tau^2=B_3^2-\tau+c,
\end{equation}
\vspace*{0.2 cm}\\
where $c$ is a constant. On the other hand, from the second equation in
~\eqref{3.3}, we
obtain $B_3^2=\frac{1}{4}-k^2-\tau^2$. Thus equation ~\eqref{eq:int}
becomes
$$
2\tau^2+\tau=C,
$$
where $C$ is a constant, and therefore also $\tau$
is constant, and we have a contradiction.
\end{proof}

From Theorem \ref{eq:theo_1} and Proposition \ref{gener-prop} we
have, in conclusion,

\begin{theorem}
\label{eq:biharmoniccurves} Let $\gamma:I\to \h_3$ be a non-geodesic
curve parametrized by arc length. Then $\gamma$ is biharmonic if
and only if
\begin{equation}\label{eq:biharmonichelix}
\left\{
\begin{array}{l}
k = constant\neq 0 ;\\
\tau = constant; \\
N_3 B_3 = 0 ;\\
k^2+\tau^2=\frac{1}{4}-B_3^2 .
\end{array}
\right.
\end{equation}
\end{theorem}
\vspace*{0.4 cm}

\section{Biharmonic helices in $\h_3$}

Now we want to determine all helices in $\h_3$ that are biharmonic
but non-geodesic. From Theorem \ref{eq:biharmoniccurves} it is
clear that, to this aim, we have to study the behaviour of $N_3$
and $B_3$.

For one thing, it follows from (\ref{eq:biharmonichelix}) that $B_3$ must be
constant. We shall show that for a curve $\gamma$ satisfying
(\ref{eq:biharmonichelix})
the constant $B_3$ cannot vanish. More precisely we prove

\begin{proposition}\label{5.1}
Let $\gamma:I\to \h_3$ be a non-geodesic curve parametrized by arc length.
If $B_3=0$, then $\tau^2=\frac{1}{4}$ and $\gamma$ is not biharmonic.
\end{proposition}

\begin{proof}
As $\gamma$ is s parametrized by arc length, we
can write
\[
T = \sin\alpha\cos\beta e_1 +
\sin\alpha\sin\beta e_2 + \cos\alpha e_3 ,
\]
where $\alpha=\alpha(s)$, $\beta=\beta(s)$. By using (\ref{eq:conn_1}) we
first have

\begin{eqnarray*}
\nabla_TT &=& (\alpha'\cos\alpha\cos\beta -  \sin\alpha
\sin\beta(\beta'-\cos\alpha)) e_1 + \\
&+& (\alpha'\cos\alpha\sin\beta +  \sin\alpha \cos\beta (\beta'-\cos\alpha))e_2
- \alpha'\sin\alpha e_3 \\
&=&kN.
\end{eqnarray*}
\vspace*{0.2 cm}\\
Next we compute $B=T\times N$, and obtain
\vspace*{0.2 cm}
\begin{eqnarray*}
B_3 = \frac{\sin^2\alpha (\beta'-\cos\alpha)}{k} .
\end{eqnarray*}
Assume now $B_3 = 0$.
We exclude the case $\sin \alpha = 0$, that implies $T = e_3$ and therefore
that
$\gamma$ is a geodesic.
Thus we must have $ \beta'-\cos\alpha = 0$, and hence
\vspace*{0.2 cm}
$$
\nabla_TT=\alpha'(\cos\alpha\cos\beta e_1+\cos\alpha\sin\beta e_2
-\sin\alpha e_3).
$$
\vspace*{0.2 cm}\\
Without loss of generality, we can assume that $\alpha' > 0$
(when $\alpha' = 0$ one has a geodesic). Then we have
\vspace*{0.2 cm}\\
\begin{align}
N &=\cos\alpha\cos\beta e_1
+\cos\alpha\sin\beta e_2-\sin\alpha e_3, \notag \\
B &= -\sin\beta e_1+\cos\beta e_2 ,\notag
\end{align}
and
\vspace*{0.2 cm}\\
$$
\nabla_TN=(- \alpha' \sin\alpha\cos\beta - \frac{1}{2}\sin\beta) e_1
+(- \alpha' \sin\alpha\sin\beta + \frac{1}{2}\cos\beta) e_2
- \alpha' \cos\alpha  e_3.
$$
\vspace*{0.2 cm}\\
Now we make use of the second Frenet equation to obtain
$$
-\tau=<\nabla_TN,B>=\frac{1}{2}
$$
\end{proof}

Thus we have
\begin{corollary}\label{coroelica}
Let $\gamma:I\to \h_3$ be a non-geodesic biharmonic
helix parametrized by arc length. Then
\begin{equation}\label{eq:biharmonichelix2}
\left\{
\begin{array}{l}
B_3 = constant\neq 0 ;\\
N_3 = 0 ;\\
k^2+\tau^2=\frac{1}{4}-B_3^2 .
\end{array}
\right.
\end{equation}
\end{corollary}
\vspace{1.5 cm}

\section{Explicit formulas for
non-geodesic biharmonic curves in $\h_3$}

In this section we use the previous results to derive the explicit parametric
equations of non-geodesic biharmonic curves in the Heisenberg group $\h_3$.\\
We first prove the following
\begin{lemma}\label{lemma1}
Let $\gamma:I\to \h_3$ be a non-geodesic curve parametrized by arc length.
If $N_3=0$, then
\begin{equation}\label{tangelica}
T(s) =\sin\alpha_0\cos\beta(s) e_1+\sin\alpha_0\sin\beta(s)
e_2+\cos\alpha_0 e_3,
\end{equation}
where $\alpha_0\in \r$.
\end{lemma}

\begin{proof}
If $\gamma'=T=T_1e_1+T_2e_2+T_3e_3$, and $\parallel T \parallel = 1$,
from
\vspace*{0.2 cm}\\
\begin{eqnarray*}
\nabla_TT &=& (T_1'+T_2T_3)e_1 + (T_2'-T_1T_3)e_2 + T_3'e_3 \\
&=& kN
\end{eqnarray*}
\vspace*{0.2 cm}\\
it follows that $N_3=0$ if and only if $T_3' = 0$, i.e. if and only if $T_3
= constant.$
Since $T_3\in [0,1]$, the lemma follows.
\end{proof}

In order to find the integral curves of the field given by
(\ref{tangelica}) that are
biharmonic but non-geodesic, we must first determine the function $\beta(s)$.
After which, a simple integration will give the wanted parametric equations.
We shall prove
\begin{theorem}
\label{eq:biharmonicintegralcurves}
The parametric equations of all non-geodesic biharmonic curves $\gamma$ of
$\h_3$ are
\begin{equation}
\label{eq:biharmonichelices1}
\left\{
\begin{array}{l}
x(s) = \frac{1}{A}\sin\alpha_0\sin(As+a) + b, \\
y(s) = -\frac{1}{A}\sin\alpha_0\cos(As+a) + c, \\
z(s) = (\cos\alpha_0+\frac{(\sin\alpha_0)^2}{2A})s  \\
\ \ \ \ \ \ \  -\frac{b}{2A}\sin\alpha_0\cos(As + a)-
\frac{c}{2A}\sin\alpha_0\sin(As + a) + d,
\end{array}
\right.
\end{equation}
 \vspace*{0.2 cm}\\
where $A =\frac{\cos\alpha_0\pm\sqrt{5(\cos\alpha_0)^2-4}}{2} $,
$\alpha_0\in(0,\arccos\frac{2\sqrt{5}}{5}]\cup[\arccos(-\frac{2\sqrt{5}}{5}),\pi)$
and $a,b,c,d \in \r$.
\end{theorem}

\begin{proof}
We shall make use of the Frenet formulas (\ref{eq:frenet}), and we shall
take into account Corollary \ref{coroelica} and Lemma \ref{lemma1}.

The covariant derivative of the vector field $T$ given by
(\ref{tangelica}) is
\vspace*{0.2 cm}\\
\begin{eqnarray*}
\nabla_TT &=&\sin\alpha_0(\cos\alpha_0-\beta')(\sin\beta e_1-\cos\beta e_2)\\
&=& kN,
\end{eqnarray*}
\vspace*{0.2 cm}\\
where $k=\vert\sin\alpha_0(\cos\alpha_0-\beta')\vert$.

Without loss of generality, we can always assume that
$\sin\alpha_0 (\cos\alpha_0-\beta') > 0$. Then we obtain
\vspace*{0.2 cm}\\
\begin{equation}
\label{eq:curburag}
k=\sin\alpha_0(\cos\alpha_0-\beta')
\end{equation}
and
\[
N=\sin\beta e_1-\cos\beta e_2.
\]
\vspace*{0.2 cm}\\
Next we have
\begin{equation}\label{eq:binorm}
B = T\times N=\cos\beta\cos\alpha_0 e_1+\sin\beta\cos\alpha_0 e_2-
\sin\alpha_0 e_3
\end{equation}
\vspace*{0.2 cm}\\
and
\vspace*{0.2 cm}
\[
\nabla_TN = \cos\beta(\beta'-\frac{1}{2}\cos\alpha_0)e_1+
\sin\beta(\beta'-\frac{1}{2}\cos\alpha_0)e_2
-\frac{1}{2}\sin\alpha_0 e_3.
\]
\vspace*{0.2 cm}\\
It follows that the geodesic torsion $\tau$ of $\gamma$ is given by
\begin{equation}
\label{eq:torsiuneag}
-\tau = <\nabla_TN,B> = (\cos\alpha_0)\beta'+\frac{1}{2}-(\cos\alpha_0)^2.
\end{equation}

If $\gamma$ is a curve with $\gamma' = T$, then this curve is non-geodesic
and biharmonic
if and only if
\begin{equation}
\label{eq:caracterizare}
\left\{
\begin{array}{l}
\beta'=constant, \\
\beta'\neq \cos\alpha_0, \\
k^2+\tau^2=\frac{1}{4}-B_3^2.
\end{array}
\right.
\end{equation}
From \eqref{eq:curburag}, \eqref{eq:binorm},
\eqref{eq:torsiuneag} and \eqref{eq:caracterizare} we obtain
\vspace*{0.2 cm}
$$
(\beta')^2-(\cos\alpha_0)\beta'+1-(\cos\alpha_0)^2=0.
$$

From this last equation we obtain
\[
\beta' = \frac{\cos\alpha_0\pm\sqrt{5(\cos\alpha_0)^2-4}}{2} = A,
\]
with the condition $(\cos\alpha_0)^2\geq\frac{4}{5}$ for reality of solutions,
 and therefore $\beta(s) = As+a$, where
 $a\in \r$.

In order to find the explicit equations for
$\gamma(s) = (x(s),y(s),z(s))$,
 we must integrate the system
 $\frac{d\gamma}{ds}=T$, that in our case is
 \vspace*{0.2 cm}
\begin{equation}
\left\{
\begin{array}{l}
\frac{dx}{ds}=\sin\alpha_0\cos(As+a),
\nonumber \\
\nonumber \\
\frac{dy}{ds}=\sin\alpha_0\sin(As+a),
\nonumber \\
\nonumber \\
\frac{dz}{ds}=\cos\alpha_0
+\frac{1}{2}\sin\alpha_0\big(\sin(As+a) x(s)-
\cos(As+a) y(s)\big ).
\end{array}
\right.
\end{equation}
\vspace*{0.2 cm}\\
The integration is immediate and yelds ~\eqref{eq:biharmonichelices1}.
\end{proof}

\begin{remark}
Biharmonic curves ~\eqref{eq:biharmonichelices1} can be obtained by
intersecting the two surfaces $S$ and $S'$ given by:
\begin{equation}
S(u,v)=
\left\{
\begin{array}{l}
x(u,v) = \frac{1}{A}\sin\alpha_0\sin(Au+a)+b, \\
y(u,v) = -\frac{1}{A}\sin\alpha_0\cos(Au+a)+c, \\
z(u,v) = v,
\end{array}
\right.
\end{equation}
and
\begin{equation}
S'(u,v)=
\left\{
\begin{array}{l}
x'(u,v) = \frac{v}{A}\sin\alpha_0\sin(Au+a)+b, \\
y'(u,v) = -\frac{v}{A}\sin\alpha_0\cos(Au+a)+c, \\
z'(u,v) = (\cos\alpha_0 + \frac{\sin^2 \alpha_0}{2 A}) u +
 \frac{b}{2} y'(u,v) -
 \frac{c}{2} x'(u,v) + d.
 \end{array}
\right.
\end{equation}
\vspace*{0.2 cm}\\
 The surface $S$ has constant non zero mean curvature; it is the ``round
cylinder''
 with rulings parallel to the axis of revolution of $\h_3$ at the point
$(b,c,0)$
 and as directrix the circle in the plane $z = 0 $ around this point;
 this circle
 has constant (geodesic) curvature also in $\h_3$. This cylinder
 has constant non zero mean curvature and zero Gaussian curvature also in
$\h_3$.
 The surface $S'$ is a ``helicoid'' which is minimal in the Heisenberg
group $\h_3$, as one
 can easily check by using the formulas given by Bekkar in \cite{MB}.
 Moreover, the biharmonic curves are geodesics of the cylinder  and the
cylinder is never a
 biharmonic surface.
 (With regard to the study of biharmonic surfaces, a paper devoted to
 invariant non-minimal biharmonic surfaces in
$\h_3$ is in preparation.)
\end{remark}

\begin{figure}[htb]
\epsfxsize=1in \centerline{
\leavevmode
\epsffile{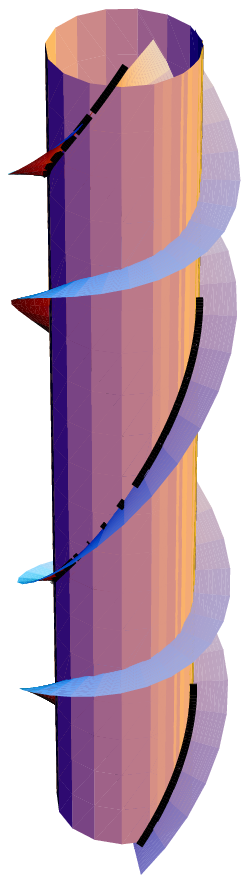}}
\protect\caption{ \label{cil-elic-biharm}}
\end{figure}

Figure 1 shows this intersection and is obtained for $a = b = c =
1$ and $\sin \alpha_0 = \frac{1}{\sqrt{10}}$.

In fact, the intersection of $S$ and $S'$ is the union of two curves.
Figure 1 shows the biharmonic one.

\begin{remark}
At each point $p\in\h_3$ the vectors tangent to biharmonic
curves form a solid cone $\mathcal{C}_p$ in $T_p\h_3$.
For each point $p \in \h_3$ and each vector $X_p \in T_p\h_3 \setminus
\mathcal{C}_p$,
the only biharmonic curve $\gamma$ arising from $p$ and such that
$\dot{\gamma}(p) = X_p$
is the geodesic determined by $p$ and $X_p$.
Thus, any $X_p \in \mathcal{C}_p$ is simultaneously tangent to a geodesic
and to a non-geodesic biharmonic curve.
We visualize this fact in the following picture, where the geodesic is the curve
not lying on the helicoid, while the interior edge is the biharmonic curve.
\end{remark}

\begin{figure}[htb]
\epsfxsize=1in
\centerline{
\leavevmode
\epsffile{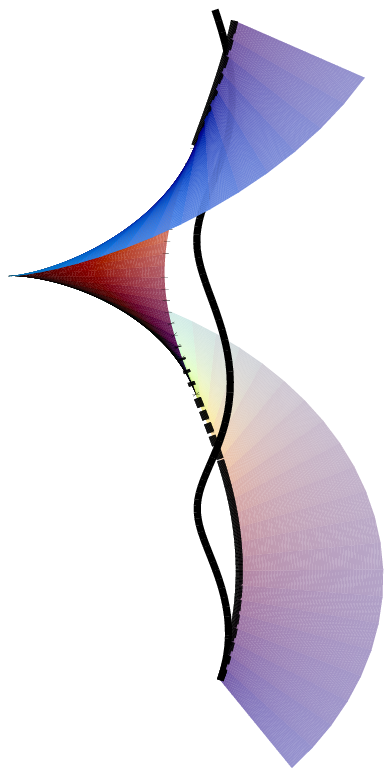}}
\protect\caption{ \label{biarm-geod}  }
\end{figure}

\begin{remark}
Let's denote by $\gamma$ the curve ~\eqref{eq:biharmonichelices1}
when $b = c = d = 0$, and by $\widetilde{\gamma}$ the
curve ~\eqref{eq:biharmonichelices1} when $b^2+c^2 > 0$.
Of course, both are helix in $\h_3$. Every biharmonic curve $\widetilde{\gamma}$
can be obtained from $\gamma$ by means of a left-translation, i.e.
$\widetilde{\gamma} = L_{(b,c,d)}\circ\gamma$, and we note
that $\gamma$ is a helix in $\r^3$, but $\widetilde{\gamma}$ is
not.
\end{remark}

\begin{remark}
The vector field $T$ tangent to a curve given by (\ref{eq:biharmonichelices1})
is transverse to the contact structure of $\h_3$  (for the contact geometry
one can
see, for example, \cite{CG}) determined by the $1$-form

\[
\theta^3 = dz  - \frac{x dy - y dx}{2}.
\]
\vspace*{0.2 cm}
(It easy to verify that the contact condition, $\theta^3 \wedge d\theta^3
\neq 0$,
is satisfied.)
In fact one has
\[
\theta^3(T) = \cos\alpha_{0} \neq 0
\]
\vspace*{0.2 cm}\\
for $\alpha_0\in(0,\arccos\frac{2\sqrt{5}}{5}]\cup
[\arccos(-\frac{2\sqrt{5}}{5}),\pi)$.
\end{remark}

Consider now a curve   $\gamma:I\to \h_3$ tangent to the contact structure
(such  a curve is called a {\it Legendre curve}), parametrized by arc
lengh. Its velocity
vector field $X$ has then the expression
\[
X = \cos\psi(s) e_{1} +  \sin\psi(s) e_{2} .
\]
\vspace*{0.2 cm}\\
It is not difficult to see that the vector fields $T$ and $X$ coincide for
 $\alpha_{0} = \frac{\pi}{2}$.
 It follows that
{\it
the Legendre curves of the Heisenberg group $H_{3}$ are biharmonic if and
only if they are
geodesic.}

\begin{remark}
The one-parameter subgroups $\sigma(u) = \exp uX$ are biharmonic if and
only if they are
geodesic. In fact, if $\sigma(u) = \exp uX$ is not a geodesic,
then $k$ and $\tau$ are always related by the formula (see \cite{PPAS})

\begin{equation}\label{4.8}
k^2 + \tau^2 = \frac{1}{4} , \hspace{2 cm} \mbox{and} \hspace{2 cm} B_3 \neq 0.
\end{equation}
Now, the assertion follows from the second equation in (\ref{3.3}).
\end{remark}

\begin{remark}
Finally we note that the methods of this paper can be extended to study
biharmonic curves
in the Cartan-Vranceanu $3$-manifolds, namely the Riemannian spaces
($\r^{3}, ds^{2}_{m,l}$),
where the Riemannian metrics  $ds^{2}_{m,l}$
are defined by
\vspace*{0.2 cm}
\begin{equation}\label{5.35}\vspace{0.3 cm}
\hspace{0.7 cm} ds^{2}_{m,l} =\frac{dx^{2} + dy^{2}}{[1 + m(x^{2} +
y^{2})]^{2}} +  \left(dz +
\frac{l}{2} \frac{ydx - xdy}{[1 + m(x^{2} + y^{2})]}\right)^{2},
\hspace{0.5 cm}l,m \in \r .
\end{equation}
\vspace*{0.2 cm}\\
This two-parameter family of metrics reduces to the metric (\ref {metricaheis})
for $m = 0 $ and $l = 1$
(see  \cite{Pi1},  \cite{Pi2},\cite{Vr} and \cite{Ca} for a discussion of
these metrics
and their properties).

The system for the non-geodesic biharmonic curves corresponding to the metric
$ds^2_{m,l}$ can be obtained by using the same techniques, and it turns out
to be

\begin{equation}\label{3.3Vran}
\left\{
\begin{array}{l}
k=constant\neq 0 ;\\
  \\
k^2+\tau^2=\frac{l^2}{4}-(l^2 - 4m) B_3^2 ;\\
 \\
\tau'= (l^2 - 4m)N_3B_3 .
\end{array}
\right.
\end{equation}
For $l^2 - 4m = 0$ the metric (\ref{5.35}) has to be of  constant curvature
$\frac{l^2}{4}$ and
we have two cases
\begin{enumerate}
\item[({\it i})] if $l = 0$, then $m= 0$, and we are in the Euclidean
space, where $\gamma(s)$ is
biharmonic if and only if it is a line;
\item[({\it ii})] if $l \neq 0$, the metric (\ref{5.35}) has constant
positive sectional curvature,
and therefore the non-geodesic biharmonic curves are circles or spherical
helices (see \cite{RCSMCO1}).
\end{enumerate}
\end{remark}

\end{document}